\documentclass[11pt,a4paper]{article}
\usepackage[a4paper,margin=27mm,headsep=6mm]{geometry}
\usepackage{amsmath,amssymb,amsthm}
\usepackage{fancyhdr}
\usepackage{tikz}
\usepackage[hidelinks]{hyperref}
\usepackage{float}
\hypersetup{
  pdftitle={A two-ended warped product manifold that does not admit the strong ratio limit property},
  pdfauthor={Bo Li},
  pdfsubject={Heat kernel asymptotics on Riemannian manifolds},
  pdfkeywords={heat kernel, strong ratio limit, warped product}
}

\newcommand{\R}{\mathbb{R}}
\newcommand{\Sphere}{\mathbb{S}}
\newcommand{\dd}{\mathrm{d}}
\newcommand{\vol}{\operatorname{vol}}

\newtheorem{theorem}{Theorem}
\newtheorem{proposition}[theorem]{Proposition}
\newtheorem{corollary}[theorem]{Corollary}
\newtheorem{lemma}[theorem]{Lemma}

\allowdisplaybreaks[2]

\title{A two-ended warped product manifold does not satisfy the strong ratio limit property}
\author{Bo Li}
\date{}

\begin{document}
\maketitle

\noindent\textbf{Abstract.}
We construct a complete two-ended warped product cylinder of bounded Gaussian curvature whose heat kernel fails the strong ratio limit property. The failure is witnessed by the quotient of two diagonal heat kernel values at fixed points. A fixed model with quadratic warping on the left and seventh-power warping on the right provides a one-sided comparison via a harmonic barrier and a Dirichlet heat kernel estimate. We then prescribe the warping on successively larger intervals. A quantitative fixed-time stability estimate allows the one-sided comparisons at the chosen times to persist in the limiting metric. The proof uses the radial heat semigroup, the parabolic maximum principle, and the Duhamel formula.


\section{Introduction and main result}

Let \(p_t(P,Q)\) be the minimal heat kernel of the Laplace--Beltrami operator on a complete connected Riemannian manifold, with respect to its Riemannian volume. Davies' strong ratio limit conjecture \cite[Section~5, (14)]{Davies} asks whether, for all fixed points \(P,Q,O\), the quotient
\[
 \frac{p_t(P,Q)}{p_t(O,O)}
\]
converges to a positive finite limit as \(t\to\infty\). In particular, the quotient of any two fixed diagonal values must converge.

Several geometric and potential-theoretic results give positive answers in specific settings. On complete manifolds with nonnegative Ricci curvature, the Li--Yau Harnack inequality \cite{LiYau} yields convergence of the normalized kernel to \(1\). For a fixed \(0<\varepsilon<1\), two applications of the inequality, one in each space variable, compare \(p_t(P,Q)\) from above to \(p_{(1+\varepsilon)t}(O,O)\) and from below to \(p_{(1-\varepsilon)t}(O,O)\). Diagonal monotonicity, followed by \(t\to\infty\) and \(\varepsilon\downarrow0\), gives the limit. Pinchover proved the ratio limit for self-adjoint elliptic operators at criticality \cite[Theorem~1.1 and Lemma~2.3]{Pinchover}. Under a time-shift bound that holds for symmetric heat kernels, he also characterized the ratio limit by uniqueness of the parabolic Martin limit along each fixed spatial pole \cite[Theorem~3.1]{Pinchover}. Ledrappier and Lim \cite[Theorem~1.1 and Corollary~1.2]{LedrappierLim} obtained the asymptotic \(p_t(P,Q)\sim c(P,Q)t^{-3/2}e^{-\lambda_0t}\) on universal covers of closed manifolds with negative sectional curvature, with \(c(P,Q)>0\) and \(\lambda_0\) the bottom of the spectrum of \(-\Delta_g\). The common time factor yields the ratio limit on these covers.

Kozma \cite[Section~2]{Kozma} disproved the graph analogue by constructing two regions that dominate the return probabilities along alternating sequences of times. He outlined a way to thicken the graph into a manifold \cite[Section~3]{Kozma}. We prove failure of the strong ratio limit property for the Laplace--Beltrami operator of a smooth two-dimensional warped product cylinder.

\begin{theorem}\label{thm:main}
There exist \(L>2\), a smooth function \(w:\R\to[1,\infty)\) with \(w=1\) on \([-2L,2L]\), and fixed observation coordinates \(x<y\) in \((-2L,2L)\) with the following property. The metric
\[
 g=\dd r^2+w(r)^2\,\dd\theta^2
\]
on \(\R\times\Sphere^1\) is complete and has bounded Gaussian curvature. Its heat kernel satisfies
\begin{equation}\label{eq:target}
 \liminf_{t\to\infty}\frac{p_t((x,0),(x,0))}{p_t((y,0),(y,0))}
 \le\frac12<2\le
 \limsup_{t\to\infty}\frac{p_t((x,0),(x,0))}{p_t((y,0),(y,0))}.
\end{equation}
In particular, the strong ratio limit property fails.
\end{theorem}

The construction realizes the alternating return mechanism through a single smooth warping function on a surface. Its flat central segment gives an explicit harmonic barrier. The value of this barrier at the opposite observation point decreases as the segment length increases. A parabolic maximum principle then controls the heat transferred between the two sides. The fixed model pairs a reflected free heat kernel on the quadratic tail with a Nash estimate on the seventh-power side. To pass from this comparison to one fixed metric, we use the Schr\"odinger form of the radial operator and a Gaussian Duhamel estimate. It allows us to change the warping farther out while controlling the heat kernel at each time already selected. A local Harnack comparison finally transfers the radial estimates to the pointwise surface kernel.

Section~\ref{sec:radial} introduces the weight and gives the radial reduction. Section~\ref{sec:model} establishes the comparison for the fixed left-quadratic, right-seventh-power model. Section~\ref{sec:locality} proves the quantitative stability estimate, and Section~\ref{sec:construction} assembles the final profile.

\section{The warped product cylinder and its radial heat kernel}\label{sec:radial}

\subsection{The weight and the metric}

Let \(L>2\). Fix a smooth nondecreasing function \(\chi:\R\to[0,1]\) equal to \(0\) on \((-\infty,0]\) and \(1\) on \([1,\infty)\). For a smooth profile \(a:\R\to[0,1]\) with \(\|a'\|_\infty\le\|\chi'\|_\infty\), define the weight \(w:\R\to[1,\infty)\) by
\begin{equation}\label{eq:weight}
\log w(r)=
\begin{cases}
\displaystyle\int_r^{-2L}\frac{\chi(-z-2L)(2+5a(z))}{1-z-2L}\,\dd z,
  &r<-2L,\\[6pt]
0,& |r|\le 2L,\\[3pt]
\displaystyle\int_{2L}^{r}\frac{\chi(z-2L)(2+5a(z))}{1+z-2L}\,\dd z,
  &r>2L.
\end{cases}
\end{equation}
On \(M=\R\times\Sphere^1\), with \(\Sphere^1=\R/(2\pi\mathbb Z)\), the metric and volume element are
\begin{equation}\label{eq:metric}
g=\dd r^2+w(r)^2\,\dd\theta^2,
\qquad \dd\vol_g=w(r)\,\dd r\,\dd\theta.
\end{equation}
We call \(M\) with this metric a two-ended warped product cylinder, in the sense of the warped product construction of Bishop and O'Neill \cite[Section~7]{BishopONeill}. Points are written \((x,\theta)\), with \(x\in\R\) and \(\theta\in\Sphere^1\). The variable \(r\) is used for the metric coefficients.

The flatness of \(\chi\) at \(0\) gives smoothness at \(r=\pm 2L\). The integrands in \eqref{eq:weight} are nonnegative. With the change of variables \(s=|z|-2L\) on each ray, we obtain
\begin{align}
 \log w(r)&\ge 2\int_1^{|r|-2L}\frac{\dd s}{1+s}
     =2\log\frac{1+|r|-2L}{2},\qquad |r|\ge 2L+1,\notag\\
 w(r)&\ge 1\quad(r\in\R),\qquad
 w(r)\ge \frac{(1+|r|-2L)^2}{4}\quad(|r|\ge 2L).\label{eq:weight-bounds}
\end{align}
For \(2L\le |r|\le 2L+1\), the second bound follows from the first because \((1+|r|-2L)^2/4\le1\). On the open rays one has
\begin{align}
(\log w)'(r)&=\operatorname{sgn}(r)\,
 \frac{\chi(|r|-2L)(2+5a(r))}{1+|r|-2L},\notag\\
\|(\log w)'\|_\infty&\le 7,\qquad
\|(\log w)''\|_\infty \le 12\|\chi'\|_\infty+7,\notag\\
K&=-\frac{w''}{w}=-(\log w)''-((\log w)')^2.\label{eq:curvature}
\end{align}
In the second derivative estimate, differentiating the cutoff, the profile and the reciprocal denominator gives bounds \(7\|\chi'\|_\infty\), \(5\|\chi'\|_\infty\) and \(7\), respectively. Thus the Gaussian curvature \(K\) is bounded uniformly over these profiles. The inequality \(w\ge1\) implies that every \(g\)-Cauchy sequence is Cauchy for the complete product metric. Its product-metric limit is also its \(g\)-limit because the two smooth metrics induce the same topology. Hence \(g\) is complete.

\subsection{Radial and pointwise heat kernels}

Write \(\mu(\dd r)=w(r)\,\dd r\). We use the nonpositive convention for the Laplace--Beltrami operator \(\Delta_g\). Its expression and its restriction to functions of \(r\) are
\begin{align}
 \Delta_g&=w^{-1}\partial_r(w\partial_r)+w^{-2}\partial_\theta^2,\notag\\
 \mathcal L f&=w^{-1}(wf')'=f''+(\log w)'f',
 \qquad
 \int_\R (-\mathcal L f)\overline f\,\dd\mu
   =\int_\R |f'|^2\,\dd\mu,\qquad f\in C_c^\infty(\R).\label{eq:radial-operator}
\end{align}
Use the closed form generated by this energy for the realization of \(-\mathcal L\) on \(L^2(\mu)\). On an open interval \(D\), use the closure of the form on \(C_c^\infty(D)\) for the Dirichlet realization \(-\mathcal L_D\). Truncation to \([0,1]\) decreases the energy, so both semigroups preserve positivity and are contractive on \(L^\infty\). By symmetry they are also contractive on \(L^1\). We write \(q_t\) and \(q_t^D\) for their respective heat kernels with respect to \(\mu\). If \(p_t\) is the heat kernel of \(e^{t\Delta_g}\) with respect to \(\vol_g\), rotation invariance and \eqref{eq:metric} give
\begin{equation}\label{eq:radial-average}
q_t(x,y)=\int_0^{2\pi}p_t((x,0),(y,\theta))\,\dd\theta .
\end{equation}
Indeed, the closed form of the surface semigroup restricted to functions of \(r\) agrees with the form in \eqref{eq:radial-operator} after normalizing the angular integral by \(2\pi\). Applying the two semigroups to the same radial function gives \eqref{eq:radial-average}.

For \(t\ge t_0>0\), symmetry and the semigroup identity give
\begin{align}
 \|q_{t_0/2}(x,\cdot)\|_{L^2(\mu)}^2 = q_{t_0}(x,x),\qquad
 q_t(x,x)=\big\langle q_{t_0/2}(x,\cdot),
 e^{(t-t_0)\mathcal L}q_{t_0/2}(x,\cdot)\big\rangle_{L^2(\mu)}.
 \label{eq:log-convex}
\end{align}
The spectral theorem for \(-\mathcal L\ge0\) represents the second line as the Laplace transform of a finite positive measure in \(t-t_0\). This representation gives monotonicity, and H\"older's inequality gives log-convexity. Since \(t_0\) is arbitrary, these properties hold for all positive times.

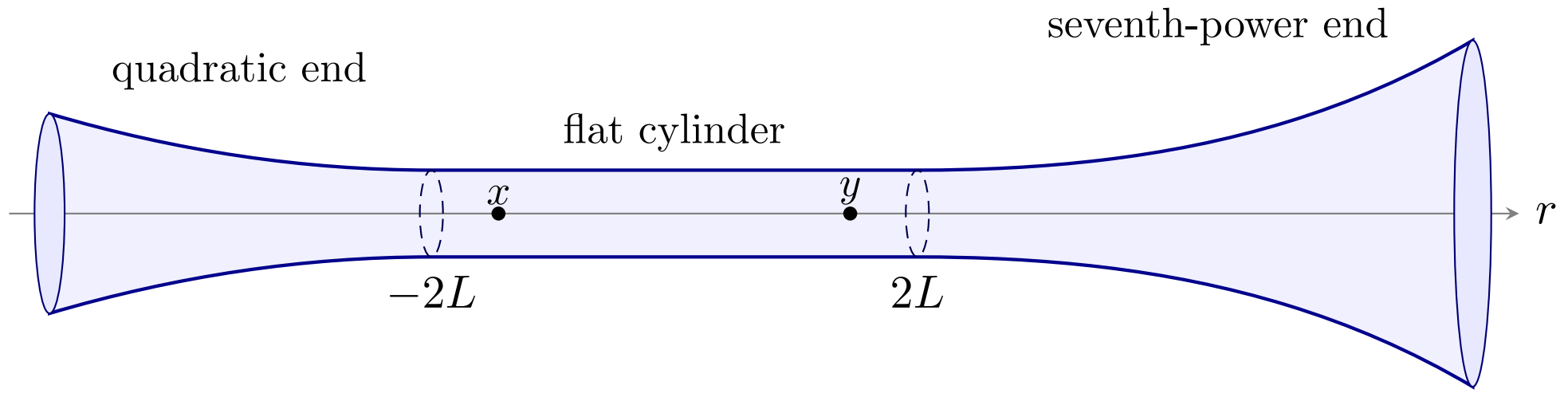
\begin{figure}[!htbp]
\centering
\begin{tikzpicture}[x=0.96cm,y=0.72cm,>=stealth]
  \filldraw[fill=blue!6,draw=blue!55!black,thick]
    (-5.4,1.15)
    .. controls (-4.3,.717) and (-3.2,.5) .. (-2.1,.5)
    -- (2.1,.5)
    .. controls (2.6333,.5) and (3.1667,.5324) .. (3.7,.6171)
    .. controls (4.2333,.7019) and (4.7667,.8389) .. (5.3,1.0572)
    .. controls (5.8333,1.2753) and (6.3667,1.5746) .. (6.9,1.995)
    -- (6.9,-1.995)
    .. controls (6.3667,-1.5746) and (5.8333,-1.2753) .. (5.3,-1.0572)
    .. controls (4.7667,-.8389) and (4.2333,-.7019) .. (3.7,-.6171)
    .. controls (3.1667,-.5324) and (2.6333,-.5) .. (2.1,-.5)
    -- (-2.1,-.5)
    .. controls (-3.2,-.5) and (-4.3,-.717) .. (-5.4,-1.15)
    -- cycle;

  \draw[->,gray] (-5.75,0)--(7.3,0) node[right,black] {\(r\)};

  \filldraw[fill=blue!9,draw=blue!45!black]
    (-5.4,0) ellipse (.13 and 1.15);
  \filldraw[fill=blue!9,draw=blue!45!black]
    (6.9,0) ellipse (.16 and 1.995);

  \draw[blue!30!black,densely dashed]
    (-2.1,0) ellipse (.10 and .50);
  \draw[blue!30!black,densely dashed]
    (2.1,0) ellipse (.10 and .50);

  \node[above] at (-3.76,1.24) {\small quadratic end};
  \node[above] at (4.7,1.75) {\small seventh-power end};
  \node[above] at (0,.53) {\small flat cylinder};
  \node[below] at (-2.1,-.53) {\(-2L\)};
  \node[below] at (2.1,-.53) {\(2L\)};

  \fill (-1.52,0) circle (1.65pt)
    node[above=2pt,inner sep=0pt,font=\small] {\(x\)};
  \fill (1.52,0) circle (1.65pt)
    node[above=2pt,inner sep=0pt,font=\small] {\(y\)};
\end{tikzpicture}
\caption{Schematic of the model cylinder with a quadratic left end
and a seventh-power right end. The dashed cross-sections mark the
boundaries of the flat region. The dots mark the axial coordinates
of the observation points \((x,0)\) and \((y,0)\).}
\label{fig:cylinder}
\end{figure}

Fix the observation coordinates \(x=-2L+2\) and \(y=2L-2\), as in Figure~\ref{fig:cylinder}.

\begin{lemma}\label{lem:pointwise}
There is an absolute constant \(C_H\ge1\), independent of the profile and of \(L\), such that, for \(z\in\{x,y\}\) and \(t\ge2\),
\begin{equation}\label{eq:pointwise}
 \frac{q_t(z,z)}{2\pi}
 \le p_t((z,0),(z,0))
 \le \frac{C_Hq_t(z,z)}{2\pi}.
\end{equation}
\end{lemma}

\begin{proof}
For the first inequality, apply the Cauchy--Schwarz inequality to the two half-time kernels and use rotation invariance of the diagonal. For the second, the positive function \(p_t((z,0),(\cdot,\cdot))\) solves the ordinary flat heat equation in the band of axial width \(4\) centered at \(z\). The local parabolic Harnack inequality \cite{Moser} gives a comparison along arcs of length at most \(\pi\) in this band over a time increment \(1\). In particular a single chain length, spatial scale and time increment per step work for all angles and both observation points. Thus
\begin{align*}
 p_t((z,0),(z,\theta))
 &\le \big[p_t((z,0),(z,0))\,p_t((z,\theta),(z,\theta))\big]^{1/2}
 =p_t((z,0),(z,0)),\\
 p_t((z,0),(z,0)) &\le C_Hp_{t+1}((z,0),(z,\theta)),\\
 q_t(z,z) &\le 2\pi p_t((z,0),(z,0))
 \le C_H q_{t+1}(z,z)\le C_H q_t(z,z).\qedhere
\end{align*}
\end{proof}

\section{Heat kernel comparison for the fixed model}\label{sec:model}

The fixed model has left terminal value \(a(r)=0\) and right terminal value \(a(r)=1\). Each condition is imposed for all \(r\) sufficiently far along the designated ray. For example, \(a(r)=\chi(r)\) has these terminal values and respects the derivative bound. On the left and right tails, respectively, \eqref{eq:weight} gives
\begin{align*}
 \frac{\dd}{\dd r}\bigl(\log w(r)-2\log(1-r-2L)\bigr)&=0,\\
 \frac{\dd}{\dd r}\bigl(\log w(r)-7\log(1+r-2L)\bigr)&=0.
\end{align*}
The weight is therefore a positive multiple of \((1-r-2L)^2\) on the left tail and of \((1+r-2L)^7\) on the right tail. We call these the quadratic end and the seventh-power end, respectively.

\subsection{The harmonic barrier}

Let \(D=(x,\infty)\). The integrability of \(w^{-1}\) on \([x,\infty)\) follows from \eqref{eq:weight-bounds}. For \(z\ge x\), set
\begin{align*}
 h(z)&=\frac{\displaystyle\int_z^\infty w(r)^{-1}\,\dd r}
              {\displaystyle\int_x^\infty w(r)^{-1}\,\dd r},&
 w(z)h'(z)&=-\left(\int_x^\infty w(r)^{-1}\,\dd r\right)^{-1},\\
 h(x)&=1,\qquad 0\le h\le1,&
 \mathcal L h&=0 \quad\text{on }D.
\end{align*}
The bound at \(y\) is uniform over profiles. The numerator includes a flat segment of length \(2\), and the whole interval \([x,y]\) contributes to the denominator. More precisely,
\begin{align}
 \int_y^\infty w(r)^{-1}\,\dd r
 &\le (2L-y)+4\int_0^\infty(1+s)^{-2}\,\dd s=2+4=6,\notag\\
 \int_x^\infty w(r)^{-1}\,\dd r
 &\ge \int_x^y 1\,\dd r=4L-4,\qquad h(y)\le\frac{6}{4L-4}.
 \label{eq:barrier-bound}
\end{align}

\subsection{The quadratic-end lower bound}

\begin{lemma}\label{lem:quadratic}
For the fixed model, there is a constant \(c>0\), depending on the profile, such that
\[
 q_t(x,x)\ge c\,t^{-3/2},\qquad t\ge3.
\]
\end{lemma}

\begin{proof}
Multiplication by \(w^{1/2}\) is unitary from \(L^2(\mu)\) to \(L^2(\dd r)\). To compute the conjugated operator, take \(f\in C_c^\infty(\R)\). The product rule gives
\begin{align}
 w(w^{-1/2}f)'&=w^{1/2}f'-(w^{1/2})'f,\notag\\
 w^{1/2}(-\mathcal L)w^{-1/2}f
 &=-w^{-1/2}\big[w^{1/2}f'-(w^{1/2})'f\big]'
 =-f''+V_wf,\notag\\
 V_w&=\frac{(w^{1/2})''}{w^{1/2}}
 =\tfrac12(\log w)''+\tfrac14((\log w)')^2.\label{eq:potential}
\end{align}
The same calculation holds on Dirichlet form domains by closure. Choose \(z<-2L\) in the quadratic tail so that \(w^{1/2}\) is affine on \((-\infty,z)\), where \(V_w=0\). Odd reflection of the free heat kernel across \(z\) and the inverse unitary transformation give
\begin{equation}\label{eq:dirichlet-reflection}
 q_t^{(-\infty,z)}(y_1,y_2)
 =\frac{e^{-(y_1-y_2)^2/(4t)}}{\sqrt{4\pi t\,w(y_1)w(y_2)}}
 \left(1-e^{-(z-y_1)(z-y_2)/t}\right),\qquad y_1,y_2<z.
\end{equation}
Here \((y_1+y_2-2z)^2-(y_1-y_2)^2=4(z-y_1)(z-y_2)\). Let \(I=[z-2,z-1]\). For \(y_1,y_2\in I\) and \(t\ge1\), the factors in \eqref{eq:dirichlet-reflection} satisfy
\begin{align*}
 e^{-(y_1-y_2)^2/(4t)}&\ge e^{-1/4},\qquad
 \sqrt{w(y_1)w(y_2)}\le\max_I w,\\
 1-e^{-(z-y_1)(z-y_2)/t}
 &=\int_0^{(z-y_1)(z-y_2)/t} e^{-s}\,\dd s
 \ge e^{-4}t^{-1}.
\end{align*}
By domain monotonicity, \(q_t(y_1,y_2)\ge q_t^{(-\infty,z)}(y_1,y_2)\ge c t^{-3/2}\) on \(I\times I\), with \(c>0\) depending on the fixed model. The semigroup identity and symmetry then give
\begin{align*}
 q_{t+2}(x,x)
 &\ge \int_I\!\!\int_I q_1(x,y_1)q_t(y_1,y_2)q_1(y_2,x)\,
      \mu(\dd y_1)\mu(\dd y_2)\\
 &\ge c t^{-3/2}\,\mu(I)^2\left(\min_{y_1\in I}q_1(x,y_1)\right)^2,
 \qquad t\ge1 .
\end{align*}
The Gaussian comparison \eqref{eq:schrodinger-comparison} below and \(w(x)=1\) give
\[
 \min_{y_1\in I}q_1(x,y_1)
 \ge \frac{\exp\big(-C_V-(x-z+2)^2/4\big)}
 {\sqrt{4\pi\max_I w}}>0.
\]
Replacing \(t\) by \(t-2\) in the semigroup lower bound proves the lemma, since \((t-2)^{-3/2}\ge t^{-3/2}\) for \(t\ge3\).
\end{proof}

\subsection{The seventh-power-end Dirichlet upper bound}

\begin{lemma}\label{lem:seventh}
Let \(D=(x,\infty)\) for the fixed model. There is a constant \(C>0\), depending on the profile, such that
\[
 q_t^D(z_1,z_2)\le C t^{-4},
 \qquad z_1,z_2\in D,\quad t>0 .
\]
\end{lemma}

\begin{proof}
There are constants \(c,C>0\) with \(c(1+r-x)^7\le w(r)\le C(1+r-x)^7\) for \(r\in D\). Indeed, the ratio between these two weights tends to a positive constant at infinity and is positive and continuous on each bounded closed interval in \([x,\infty)\). For \(f\in C_c^\infty(D)\) take the radial function \(F(\xi)=f(x+|\xi|-1)\) for \(|\xi|>1\) in \(\R^8\), and \(F=0\) for \(|\xi|\le1\). With \(|\Sphere^7|\) denoting the unit sphere area, polar integration gives
\begin{align*}
 \|F\|_{L^p(\R^8)}^p &=|\Sphere^7|\int_x^\infty |f(r)|^p(1+r-x)^7\,\dd r,
 \quad p=1,2,\\
 \|\nabla F\|_{L^2(\R^8)}^2 &=|\Sphere^7|\int_x^\infty |f'(r)|^2(1+r-x)^7\,\dd r .
\end{align*}
The Euclidean Nash inequality \cite{Nash}, with exponent \(2+4/8=5/2\) on the \(L^2\) norm, now reads
\begin{equation}\label{eq:nash}
 \|f\|_{L^2(D,\mu)}^{5/2}
 \le C\,\|f\|_{L^1(D,\mu)}^{1/2}\int_D|f'|^2\,\dd\mu.
\end{equation}
The inequality extends to the Dirichlet form domain in \(L^1(D,\mu)\) by truncation and compactly supported form approximation. For \(f\in L^1\cap L^2(D,\mu)\), \(f\ne0\), set \(Y(s)=\|e^{s\mathcal L_D}f\|_2^2\), where the norms refer to \(\mu\) on \(D\). The form identity, \eqref{eq:nash}, and \(L^1\)-contractivity give
\begin{align*}
 Y'(s)&=-2\int_D |(e^{s\mathcal L_D}f)'|^2\,\dd\mu
 \le -c\|f\|_1^{-1/2}Y(s)^{5/4},\\
 \frac{\dd}{\dd s}\big(Y(s)^{-1/4}\big)
 &\ge c\|f\|_1^{-1/2},\qquad
 Y(s)\le C s^{-4}\|f\|_1^2.
\end{align*}
By duality and the semigroup identity, the \(L^1\)-to-\(L^2\) estimate \(\|e^{s\mathcal L_D}\|_{1\to2}\le Cs^{-2}\) yields
\[
 \|e^{t\mathcal L_D}\|_{1\to\infty}
 \le \|e^{(t/2)\mathcal L_D}\|_{2\to\infty}\,
       \|e^{(t/2)\mathcal L_D}\|_{1\to2}
 \le C t^{-4}.
\]
The stated pointwise bound follows from the operator norm by approximation of point masses and continuity of the kernel.
\end{proof}

\subsection{Comparison across the flat cylinder}

\begin{proposition}\label{prop:model}
For the fixed model with the quadratic end on the left and the seventh-power end on the right, its radial heat kernel satisfies
\begin{equation}\label{eq:model-comparison}
 \limsup_{t\to\infty}\frac{q_t(y,y)}{q_t(x,x)}\le
 \left(\frac{6}{4L-4}\right)^2 .
\end{equation}
\end{proposition}

\begin{proof}
Fix \(t>4\) and a nonnegative \(f\in C_c^\infty(\R)\). For \(0\le s\le\sqrt t\) and \(z\in D\), define
\begin{align*}
 u(s,z)&=(e^{(t/2-\sqrt t+s)\mathcal L}f)(z),\\
 v(s,z)&=\big(e^{s\mathcal L_D}[(e^{(t/2-\sqrt t)\mathcal L}f)|_D]\big)(z),\\
 u(s,x)&\le\|q_{t/2-\sqrt t+s}(x,\cdot)\|_{L^2(\mu)}\|f\|_{L^2(\mu)}
 \le \sqrt{q_{t-2\sqrt t}(x,x)}\|f\|_{L^2(\mu)}=B_f .
\end{align*}
The functions \(u\) and \(v\) solve the same heat equation with the same initial value on \(D\), while \(v(s,x)=0\) for \(s>0\). The last line follows from \eqref{eq:log-convex} and diagonal monotonicity. The harmonic function \(B_fh\) controls \(u-v\) at \(x\), including the upper limit at the corner \((0,x)\) since \(v\ge0\). To impose the comparison at infinity, for \(\varepsilon>0\) subtract \(\varepsilon e^{10s}(1+z^2)\) from \(u-v-B_fh\). By \eqref{eq:curvature},
\begin{align*}
 (\partial_s-\mathcal L)\big(e^{10s}(1+z^2)\big)
 &=e^{10s}(8+10z^2-2z(\log w)'(z))\\
 &\ge e^{10s}\big(10(|z|-7/10)^2+31/10\big)>0 .
\end{align*}
Since \(u\le\|f\|_\infty\), \(v\ge0\), and \(h\ge0\), the subtracted function is nonpositive at any sufficiently large right endpoint, uniformly in \(s\in[0,\sqrt t]\). The maximum principle on each such finite interval, followed by \(\varepsilon\downarrow0\), gives \(u-v\le B_fh\) throughout the half-line.

At \(s=\sqrt t\) and \(z=y\), use the Cauchy--Schwarz inequality and the \(L^2\)-contractivity of \(e^{(t/2-\sqrt t)\mathcal L}\) to obtain
\begin{align}
 v(\sqrt t,y)
 &\le \|q_{\sqrt t}^D(y,\cdot)\|_{L^2(D,\mu)}
       \big\|(e^{(t/2-\sqrt t)\mathcal L}f)|_D\big\|_{L^2(D,\mu)}
 \le \sqrt{q_{2\sqrt t}^D(y,y)}\|f\|_{L^2(\mu)},\notag\\
 (e^{(t/2)\mathcal L}f)(y)
 &\le\left(\sqrt{q_{2\sqrt t}^D(y,y)}
       +h(y)\sqrt{q_{t-2\sqrt t}(x,x)}\right)\|f\|_{L^2(\mu)},\notag\\
 \sqrt{q_t(y,y)}
 &\le \sqrt{q_{2\sqrt t}^D(y,y)}+h(y)\sqrt{q_{t-2\sqrt t}(x,x)} .
 \label{eq:decomposition}
\end{align}
In the last line we take the supremum over the nonnegative test functions of \(L^2(\mu)\)-norm at most \(1\). By density this supremum equals \(\|q_{t/2}(y,\cdot)\|_{L^2(\mu)}=\sqrt{q_t(y,y)}\).

Using Lemmas~\ref{lem:quadratic} and~\ref{lem:seventh} in the first line and log-convexity on \([t/2,t]\) in the second line gives
\begin{align*}
 0\le\frac{q_{2\sqrt t}^D(y,y)}{q_t(x,x)}
 &\le \frac{C(2\sqrt t)^{-4}}{ct^{-3/2}}\le C t^{-1/2},\\
 0\le\log\frac{q_{t-2\sqrt t}(x,x)}{q_t(x,x)}
 &\le \frac{4}{\sqrt t}\log\frac{q_{t/2}(x,x)}{q_t(x,x)}
 \le \frac{4}{\sqrt t}\log\frac{q_1(x,x)}{ct^{-3/2}}\longrightarrow0 .
\end{align*}
Here \(t\) is sufficiently large that \(t/2\le t-2\sqrt t\le t\). Dividing \eqref{eq:decomposition} by \(\sqrt{q_t(x,x)}\), taking the upper limit and applying \eqref{eq:barrier-bound} prove \eqref{eq:model-comparison}.
\end{proof}

\section{Stability under remote changes of the profile}\label{sec:locality}

Let \(H_w=-\partial_r^2+V_w\), with \(V_w\) given by \eqref{eq:potential}, and denote its heat kernel with respect to Lebesgue measure by \(k_t^w\). Equations \eqref{eq:curvature} and \eqref{eq:potential} give the same finite bound \(C_V>0\) for \(\|V_w\|_\infty\) for every admissible profile. Write
\[
 G_t(x,y)=\frac{1}{\sqrt{4\pi t}}\exp\!\left(-\frac{(x-y)^2}{4t}\right)
\]
for the free heat kernel. The inverse unitary transformation gives the first line below. For nonnegative functions the multiplication operator \(e^{-tV_w/n}\) lies between multiplication by \(e^{-tC_V/n}\) and multiplication by \(e^{tC_V/n}\). Iterating these inequalities with the positive free semigroup and taking the strong \(L^2\) limit by the Trotter product formula gives the second line, first as an operator comparison on nonnegative functions and then pointwise by continuity.
\begin{align}
 q_t(x,y)&=\frac{k_t^w(x,y)}{\sqrt{w(x)w(y)}},\notag\\
 e^{-C_Vt}G_t(x,y)&\le k_t^w(x,y)\le e^{C_Vt}G_t(x,y).
 \label{eq:schrodinger-comparison}
\end{align}

\begin{lemma}\label{lem:locality}
Suppose the weights \(w,\widetilde w\) of two admissible profiles agree on \([-\rho,\rho]\), with \(\rho>2L\). Let \(q_t,\widetilde q_t\) be their radial kernels. For \(z\in\{x,y\}\) and \(t>0\),
\begin{equation}\label{eq:locality}
 |q_t(z,z)-\widetilde q_t(z,z)|
 \le \frac{2\sqrt{2}\,C_Vt\,e^{C_Vt}}{\sqrt{4\pi t}}
       \exp\!\left(-\frac{(\rho-|z|)^2}{2t}\right).
\end{equation}
\end{lemma}

\begin{proof}
The potentials agree on \([-\rho,\rho]\). Differentiating \(e^{-(t-s)H_w}e^{-sH_{\widetilde w}}\) and integrating in \(s\) yields the bounded-perturbation formula \cite[Section~1]{BDS}. The corresponding kernel integral is absolutely convergent by \eqref{eq:schrodinger-comparison} and Gaussian convolution. At each \(z\in\{x,y\}\), both weights have value \(1\), hence
\begin{align}
 k_t^w(z,z)-k_t^{\widetilde w}(z,z)
 &=\int_0^t\!\int_{\R} k_{t-s}^w(z,r)
  (V_{\widetilde w}-V_w)(r)\ k_s^{\widetilde w}(r,z)\,\dd r\,\dd s,\notag\\
 |q_t(z,z)-\widetilde q_t(z,z)|
 &\le 2C_V e^{C_Vt}\int_0^t\!\int_{|r|\ge\rho}
          G_{t-s}(z,r)G_s(r,z)\,\dd r\,\dd s.
 \label{eq:duhamel}
\end{align}
For \(0<s<t\), the identity for the product of Gaussians and \(s(t-s)\le t^2/4\) give, when \(|r|\ge\rho\),
\begin{align*}
 G_{t-s}(z,r)G_s(r,z)
 &=\frac{\exp\!\left(-t(r-z)^2/[4s(t-s)]\right)}{4\pi\sqrt{s(t-s)}}\\
 &\le e^{-(\rho-|z|)^2/(2t)}
   \frac{\exp\!\left(-t(r-z)^2/[8s(t-s)]\right)}{4\pi\sqrt{s(t-s)}} .
\end{align*}
The integral of the last Gaussian with its prefactor over \(\R\) equals
\[
 \frac{\sqrt{8s(t-s)/t}}{4\pi\sqrt{s(t-s)}}
      \int_{\R}e^{-u^2}\,\dd u
 =\frac{\sqrt{2}}{\sqrt{4\pi t}}.
\]
Substitution into \eqref{eq:duhamel} proves \eqref{eq:locality}. In particular the bound tends to zero as \(\rho\to\infty\) for fixed \(t\), uniformly over the two profiles being compared.
\end{proof}

\section{Construction of the counterexample}\label{sec:construction}

\begin{figure}[tbp]
\centering
\begin{tikzpicture}[x=0.93cm,y=0.72cm,quadratic/.style={fill=blue!12},seventh/.style={fill=orange!19}]
 \draw[quadratic] (-5.4,-0.29) rectangle (-4.15,0.29);
 \draw[seventh] (-4.15,-0.29) rectangle (-2.9,0.29);
 \draw[quadratic] (-2.9,-0.29) rectangle (-1.65,0.29);
 \draw[fill=gray!8] (-1.65,-0.29) rectangle (1.65,0.29);
 \draw[seventh] (1.65,-0.29) rectangle (2.9,0.29);
 \draw[quadratic] (2.9,-0.29) rectangle (4.15,0.29);
 \draw[seventh] (4.15,-0.29) rectangle (5.4,0.29);
 \node at (-5.82,0) {\(\cdots\)};
 \node at (5.82,0) {\(\cdots\)};
 \node at (-4.78,0) {\(2\)};
 \node at (-3.53,0) {\(7\)};
 \node at (-2.28,0) {\(2\)};
 \node at (0,0) {\(w=1\)};
 \node at (2.28,0) {\(7\)};
 \node at (3.53,0) {\(2\)};
 \node at (4.78,0) {\(7\)};
 \node[below] at (-4.15,-0.30) {\(-\rho_2\)};
 \node[below] at (-2.9,-0.30) {\(-\rho_1\)};
 \node[below] at (-1.65,-0.30) {\(-2L\)};
 \node[below] at (1.65,-0.30) {\(2L\)};
 \node[below] at (2.9,-0.30) {\(\rho_1\)};
 \node[below] at (4.15,-0.30) {\(\rho_2\)};
\end{tikzpicture}
\caption{Power coefficients \(2+5a(r)\) on the successive outer intervals of the final cylinder. The vertical marks locate the transition intervals schematically.}\label{fig:profile}
\end{figure}
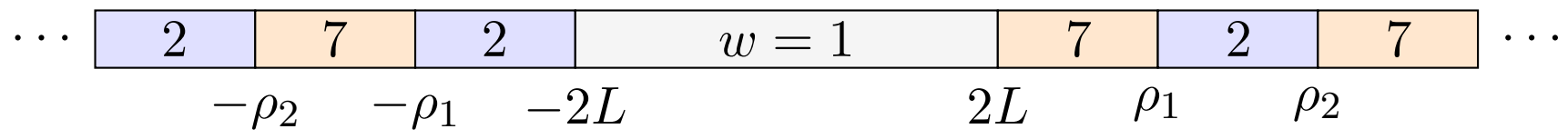

\begin{proof}[Proof of Theorem~\ref{thm:main}]
Choose \(L>2\) large enough that
\[
 \left(\frac{6}{4L-4}\right)^2<\frac{1}{32C_H},
 \qquad \delta=\frac{1}{16C_H}.
\]
Here \(C_H\) is the absolute constant in Lemma~\ref{lem:pointwise}. The observation coordinates \(x=-2L+2\) and \(y=2L-2\) stay fixed throughout the construction.

We prescribe one profile \(a_n:\R\to[0,1]\) at each stage and let \(w_n\) and \(q_t^{(n)}\) be its weight and radial kernel. For odd \(n\) the terminal values of \(a_n\) are \(0\) on the left and \(1\) on the right. For even \(n\) the terminal values are \(1\) on the left and \(0\) on the right. The reflected auxiliary weight \(w_n(-r)\) then has the orientation of the fixed model. Set
\[
 (x_n,y_n)=
 \begin{cases}
 (x,y),& n\ \text{odd},\\
 (y,x),& n\ \text{even}.
 \end{cases}
\]
Reflection preserves the divergence form, with the reflected measure, and interchanges \(x\) and \(y\). Thus Proposition~\ref{prop:model} gives, for every \(n\),
\begin{equation}\label{eq:stage-limit}
 \limsup_{t\to\infty}\frac{q_t^{(n)}(y_n,y_n)}{q_t^{(n)}(x_n,x_n)}
 \le \left(\frac{6}{4L-4}\right)^2<\frac{\delta}{2}.
\end{equation}

Start with \(a_1(r)=\chi(r)\). Once \(a_n\) is known, choose \(t_n\ge\max\{2,n\}\), strictly larger than the previous times, by \eqref{eq:stage-limit}. Next choose \(\rho_n>2L\), beyond the start of both constant tails, with \(\rho_n>\rho_{n-1}+2\) if \(n\ge2\). Increase \(\rho_n\) if necessary using Lemma~\ref{lem:locality}. The resulting choices satisfy
\begin{align}
 q_{t_n}^{(n)}(y_n,y_n)&\le\delta\,q_{t_n}^{(n)}(x_n,x_n),\notag\\
 \big|\widetilde q_{t_n}(z,z)-q_{t_n}^{(n)}(z,z)\big|
 &\le\delta\,q_{t_n}^{(n)}(x_n,x_n),\qquad z\in\{x,y\},
 \label{eq:stage}
\end{align}
where the second inequality holds for every admissible weight \(\widetilde w\) that agrees with \(w_n\) on \([-\rho_n,\rho_n]\). This choice of \(\rho_n\) is possible since \(q_{t_n}^{(n)}(x_n,x_n)>0\) and the right-hand side of \eqref{eq:locality} tends to zero for fixed \(t_n\).

Let \(\alpha_n=(1+(-1)^n)/2\), the terminal value of \(a_n\) on the left. Define the next profile by
\[
 a_{n+1}(r)=
 \begin{cases}
 \alpha_n+(\alpha_{n+1}-\alpha_n)\chi(-r-\rho_n),& r<-\rho_n,\\
 a_n(r),& |r|\le\rho_n,\\
 1-\alpha_n+(\alpha_n-\alpha_{n+1})\chi(r-\rho_n),& r>\rho_n.
 \end{cases}
\]
The two outer pieces change to their new terminal values over intervals of length \(1\). Since \(a_n\) is constant near both joining points and \(\chi\) is flat at \(0\) and \(1\), the function \(a_{n+1}\) is smooth and obeys \(\|a_{n+1}'\|_\infty\le\|\chi'\|_\infty\). By \eqref{eq:weight}, \(w_{n+1}=w_n\) on \([-\rho_n,\rho_n]\). Thus the profiles stabilize on every compact interval, defining a smooth admissible \(a\) and an associated \(w\). Figure~\ref{fig:profile} shows the resulting sequence of power coefficients \(2+5a(r)\) along the two rays.

For the radial kernel \(q_t\) of \(w\), the agreement \(w=w_n\) on \([-\rho_n,\rho_n]\) allows us to use \eqref{eq:stage} at time \(t_n\). By \eqref{eq:pointwise} and \(\delta\le1/16\),
\begin{align*}
 \frac{q_{t_n}(y_n,y_n)}{q_{t_n}(x_n,x_n)}
 &\le \frac{2\delta}{1-\delta}\le\frac{2}{15C_H},\\
 \frac{p_{t_n}((y_n,0),(y_n,0))}{p_{t_n}((x_n,0),(x_n,0))}
 &\le C_H\,\frac{q_{t_n}(y_n,y_n)}{q_{t_n}(x_n,x_n)}
 \le \frac{2}{15}.
\end{align*}
For even \(n\), the quotient in \eqref{eq:target} is at most \(2/15<1/2\). For odd \(n\), inversion gives a lower bound of \(15/2>2\). Since \(t_n\to\infty\), these two bounds prove \eqref{eq:target}. Completeness and bounded Gaussian curvature follow from the uniform estimates in Section~\ref{sec:radial}.
\end{proof}

The strong ratio limit property can fail even on a manifold with
only one topological end.

\begin{corollary}
There is a complete Riemannian three-manifold with one topological end
whose heat kernel fails the strong ratio limit property.
\end{corollary}

\begin{proof}
Let \(M=\mathbb R\times\mathbb S^1\) be the two-ended warped product
cylinder in Theorem~\ref{thm:main}, and let \(x<y\) be the two
observation coordinates in that theorem. Take
\[
 N=\mathbb R\times M=\mathbb R^2\times\mathbb S^1,
 \qquad
 g_N=\mathrm du^2+g
     =\mathrm du^2+\mathrm dr^2+w(r)^2\mathrm d\theta^2,
\]
where \(u\) is the coordinate on the first factor. We write points in
the order \((u,r,\theta)\). The product is three-dimensional and
complete because both factors are complete.

For each integer \(n\ge1\), let \(K_n=[-n,n]^2\times\mathbb S^1\).
Every compact subset of \(N\) is contained in some \(K_n\), and
\[
 N\setminus K_n
 =
 \bigl(\mathbb R^2\setminus[-n,n]^2\bigr)\times\mathbb S^1.
\]
The complement of the square is the union of the four open
half-planes \(r<-n\), \(u>n\), \(r>n\), and \(u<-n\). Each
half-plane is path-connected and intersects the next one in this
cyclic order. Their union is therefore path-connected.

Take any two points \((u_1,r_1,\theta_1)\) and
\((u_2,r_2,\theta_2)\) in \(N\setminus K_n\). First keep the angle
equal to \(\theta_1\) and move from \((u_1,r_1)\) to \((u_2,r_2)\)
along a path outside the square. This joins
\((u_1,r_1,\theta_1)\) to \((u_2,r_2,\theta_1)\). Next keep
\((u_2,r_2)\) fixed and move along the circle from \(\theta_1\) to
\(\theta_2\). The entire circle over \((u_2,r_2)\) lies in the
complement, so the second step also stays in \(N\setminus K_n\).
Figure~\ref{fig:one-topological-end} illustrates these two steps.

The complement is unbounded, since it contains a ray with \(u=0\)
and fixed \(\theta\) along which \(r\to+\infty\), and lengths in
\(N\) dominate changes in \(r\). Thus each complement in this
compact exhaustion has a single unbounded connected component.
It follows that \(N\) has one topological end.

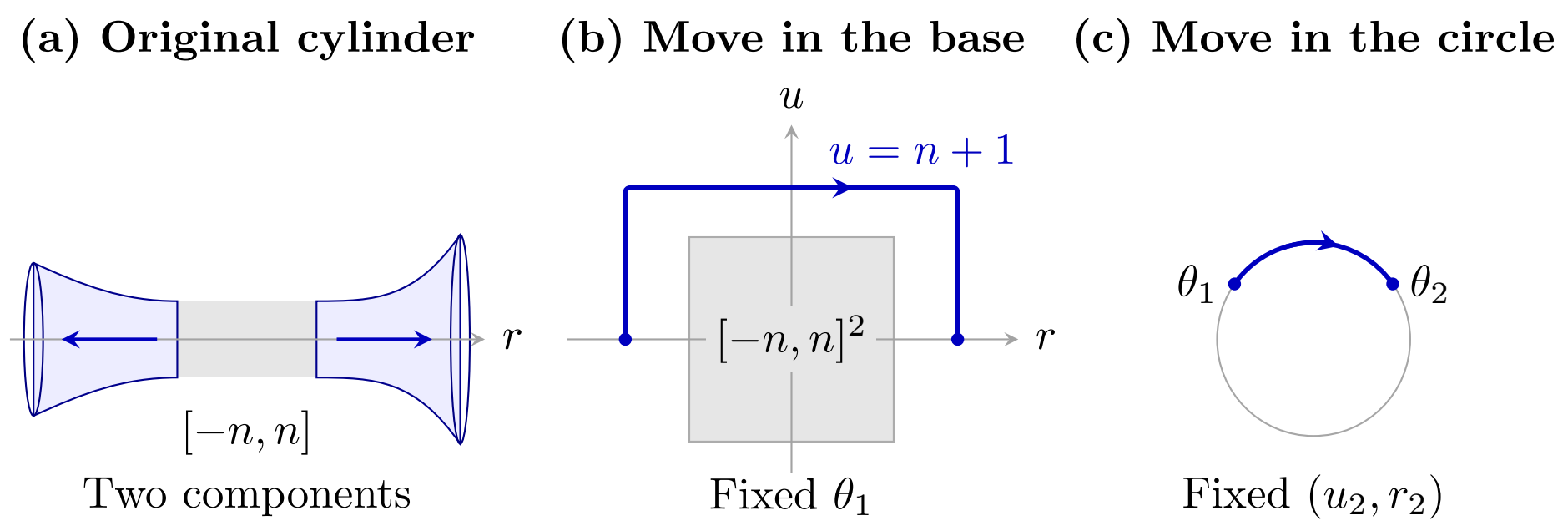
\begin{figure}[H]
\centering
\begin{tikzpicture}[x=.89cm,y=.89cm,>=stealth,font=\small]
  \begin{scope}[shift={(-4.85,0)}]
    \node[font=\small\bfseries] at (0,2.66)
      {(a) Original cylinder};

    \fill[gray!20] (-.62,-.34) rectangle (.62,.34);
    \filldraw[fill=blue!7,draw=blue!55!black]
      (-1.90,.68)
      .. controls (-1.46,.48) and (-1.10,.34) .. (-.62,.34)
      -- (-.62,-.34)
      .. controls (-1.10,-.34) and (-1.46,-.48) .. (-1.90,-.68)
      -- cycle;
    \filldraw[fill=blue!7,draw=blue!55!black]
      (.62,.34)
      .. controls (1.15,.34) and (1.50,.36) .. (1.90,.93)
      -- (1.90,-.93)
      .. controls (1.50,-.36) and (1.15,-.34) .. (.62,-.34)
      -- cycle;
    \draw[blue!50!black] (-1.90,0) ellipse (.085 and .68);
    \draw[blue!50!black] (1.90,0) ellipse (.085 and .93);
    \draw[gray!72,->] (-2.11,0)--(2.12,0)
      node[right,black] {\(r\)};
    \draw[blue!75!black,thick,->] (-.80,0)--(-1.65,0);
    \draw[blue!75!black,thick,->] (.80,0)--(1.65,0);
    \node[below] at (0,-.48) {\([-n,n]\)};
    \node[align=center] at (0,-1.40) {Two components};
  \end{scope}

  \begin{scope}
    \node[font=\small\bfseries] at (0,2.66)
      {(b) Move in the base};

    \filldraw[fill=gray!20,draw=gray!65]
      (-.91,-.91) rectangle (.91,.91);
    \draw[gray!72,->] (-2,0)--(2.02,0)
      node[right,black] {\(r\)};
    \draw[gray!72,->] (0,-1.19)--(0,1.91)
      node[above,black] {\(u\)};
    \node[fill=gray!20,inner sep=2pt] at (0,0)
      {\([-n,n]^2\)};
    \draw[blue!75!black,line width=1.05pt,rounded corners=1pt]
      (-1.48,0)--(-1.48,1.35)--(1.48,1.35)--(1.48,0);
    \draw[blue!75!black,line width=1.05pt,->]
      (-.62,1.35)--(.55,1.35);
    \fill[blue!75!black] (-1.48,0) circle (1.45pt);
    \fill[blue!75!black] (1.48,0) circle (1.45pt);
    \node[above right,text=blue!75!black] at (.18,1.38)
      {\(u=n+1\)};
    \node[align=center] at (0,-1.40) {Fixed \(\theta_1\)};
  \end{scope}

  \begin{scope}[shift={(4.65,0)}]
    \node[font=\small\bfseries] at (0,2.66)
      {(c) Move in the circle};

    \draw[gray!70] (0,0) circle (.86);
    \draw[blue!75!black,line width=1.05pt]
      (145:.86) arc (145:35:.86);
    \draw[blue!75!black,line width=1.05pt,->]
      (108:.86) arc (108:75:.86);
    \fill[blue!75!black] (145:.86) circle (1.45pt)
      node[left,black] {\(\theta_1\)};
    \fill[blue!75!black] (35:.86) circle (1.45pt)
      node[right,black] {\(\theta_2\)};
    \node[align=center] at (0,-1.40) {Fixed \((u_2,r_2)\)};
  \end{scope}
\end{tikzpicture}
\caption{Panel (a) shows the two components of the original cylinder
after removing the middle compact segment. Panel (b) shows a path
around the square in the product, with the angle held fixed.
Panel (c) shows how to change the angle once the base point has
reached its destination.}
\label{fig:one-topological-end}
\end{figure}

The Laplace--Beltrami operator on \(N\) is
\[
 \Delta_N=\partial_u^2+\Delta_g.
\]
The summands act on separate factors, so the heat kernel on \(N\)
is the product of the heat kernels on \(\mathbb R\) and \(M\).
Writing \(p_t^N\) and \(p_t\) for the kernels on \(N\) and \(M\),
respectively, we obtain for every \(t>0\)
\[
\frac{p_t^N((0,x,0),(0,x,0))}
      {p_t^N((0,y,0),(0,y,0))}
= \frac{(4\pi t)^{-1/2}p_t((x,0),(x,0))}
      {(4\pi t)^{-1/2}p_t((y,0),(y,0))}
= \frac{p_t((x,0),(x,0))}{p_t((y,0),(y,0))}.
\]
The quotient on the right does not converge as \(t\to\infty\)
by Theorem~\ref{thm:main}. Hence the strong ratio limit property
fails on \(N\).
\end{proof}

\medskip
\noindent\textbf{Acknowledgments.} 
The author acknowledges the use of ChatGPT in the preparation of this manuscript and assumes full responsibility for its content.
Bo Li is supported by the National Natural Science Foundation of China (Nos.~12471094 and 12571121) and the Natural Science Foundation of Zhejiang Province (No.~LMS26A010015).

\medskip
\noindent Bo Li (\href{mailto:bli@zjxu.edu.cn}{bli@zjxu.edu.cn}), Institute of Mathematics, Jiaxing University, Jiaxing 314001, China.

\end{document}